\documentclass[reqno]{amsart}
\usepackage{amsmath,amssymb}

\theoremstyle{plain}
\swapnumbers

\numberwithin{equation}{section}

\def\cal{\mathcal}

\newtheorem{theorem}{Theorem}[section]
\newtheorem{corollary}[theorem]{Corollary}

\newtheorem{lemma}[theorem]{Lemma}

\def\N#1#2{{\cal N}_{#1,#2}}        
\def\dualN#1#2{{\cal N}'_{#1,#2}}   
\def\Nom#1{{\cal N}_{#1}}           
\def\dualNom#1{{\cal N}'_{#1}}      
\def\XDX#1#2#3{X_{#1} \Delta_{#2} X_{#3}}       
\def\DXD#1#2#3{\Delta_{#1} X_{#2} \Delta_{#3}}  
\def\BM{{\cal B}}                   
\def\cV{{\cal V}}                   
\def\M{{\cal M}}                    
\def\0{{\bf 0}}                     
\def\Sinv#1{{#1}^{(-)}}             
\def\SinvT#1{{#1}^{(-)T}}           
\def\inv#1{{#1}^{-1}}               
\def\invT#1{({#1}^{-1})^T}              
\def\schur{\circ}
\def\com{\mathbb{C}}                
\def\D{\Delta} 
\def\T#1#2{\Theta_{#1,#2}}          
\def\Evector#1#2#3#4{#1 e_{#2} \schur #3 e_{#4}} 
\def\range#1#2#3{#1=#2, \ldots, #3}  

\def\mat#1#2{{\mathbf{M}_{#1}(#2)}} 
\def\nbyn{n \times n}               

\def\1{{\mathbf{1}}}                 
\def\2{{\mathbf{2}}}                 
\def\3{{\mathbf{3}}}                 
\def\4{{\mathbf{4}}}                 
\def\y#1#2{\mathbf{Y}^{#1,#2}_{i,j}}   

\def\Evone{
\left(
\begin{aligned}
\Evector{A}{i}{\Sinv{A}}{j}\\ 
\Evector{A}{i}{\Sinv{A}}{j}\\ 
\Evector{\SinvT{B}}{i}{B^T}{j} \\ 
\Evector{\SinvT{B}}{i}{B^T}{j} \\
\end{aligned}
\right)
}    

\def\Evtwo{
\left(
\begin{aligned}
-&\Evector{A}{i}{\Sinv{A}}{j}\\ 
-&\Evector{A}{i}{\Sinv{A}}{j}\\ 
&\Evector{\SinvT{B}}{i}{B^T}{j} \\ 
&\Evector{\SinvT{B}}{i}{B^T}{j} \\
\end{aligned}
\right)
}    

\def\Evthree{
\left(
\begin{aligned}
\Evector{\Sinv{B}}{i}{B}{j} \\ 
\Evector{\Sinv{B}}{i}{B}{j} \\
\Evector{A}{i}{\Sinv{A}}{j}\\ 
\Evector{A}{i}{\Sinv{A}}{j}\\ 
\end{aligned}
\right)
}    

\def\Evfour{
\left(
\begin{aligned}
&\Evector{\Sinv{B}}{i}{B}{j} \\ 
&\Evector{\Sinv{B}}{i}{B}{j} \\
-&\Evector{A}{i}{\Sinv{A}}{j}\\ 
-&\Evector{A}{i}{\Sinv{A}}{j}\\ 
\end{aligned}
\right)
}    

\def\Evfive{
\left(
\begin{aligned}
d & \Evector{A}{i}{B}{j} \\ 
-d & \Evector{A}{i}{B}{j} \\
d^{-1} & \Evector{\SinvT{B}}{i}{\Sinv{A}}{j} \\ 
-d^{-1} & \Evector{\SinvT{B}}{i}{\Sinv{A}}{j} \\
\end{aligned}
\right)
}    

\def\Evsix{
\left(
\begin{aligned}
d & \Evector{A}{i}{B}{j} \\ 
-d & \Evector{A}{i}{B}{j} \\
-d^{-1} & \Evector{\SinvT{B}}{i}{\Sinv{A}}{j} \\ 
d^{-1} & \Evector{\SinvT{B}}{i}{\Sinv{A}}{j} \\
\end{aligned}
\right)
}    

\def\Evseven{
\left(
\begin{aligned}
d^{-1} & \Evector{\Sinv{B}}{i}{\Sinv{A}}{j} \\ 
-d^{-1} & \Evector{\Sinv{B}}{i}{\Sinv{A}}{j} \\
d & \Evector{A}{i}{B^T}{j} \\ 
-d & \Evector{A}{i}{B^T}{j} \\
\end{aligned}
\right)
}    

\def\Eveight{
\left(
\begin{aligned}
d^{-1} & \Evector{\Sinv{B}}{i}{\Sinv{A}}{j} \\ 
-d^{-1} & \Evector{\Sinv{B}}{i}{\Sinv{A}}{j} \\
-d & \Evector{A}{i}{B^T}{j} \\ 
d & \Evector{A}{i}{B^T}{j} \\
\end{aligned}
\right)
}    

\def\sqr#1#2{{\vbox{\hrule height.#2pt
    \hbox{\vrule width.#2pt height#1pt \kern#1pt
        \vrule width.#2pt}\hrule height.#2pt}}}
\def\eqed{\sqr73}
\def\qed{%
    \ifmmode\eqno\eqed
    \else\nobreak\ \hfill\eqed\medbreak\fi}

\begin{document}
\title{Bose-Mesner Algebras attached to Invertible Jones Pairs}

\author{Ada Chan}
\address{Department of Combinatorics and Optimization,
University of Waterloo, Waterloo, Ontario, Canada N2L 3G1}
\curraddr{Department of Mathematics, 
California Institute of Technology, Pasadena, California 91125, U.S.A.}
\email{ssachan@alumni.uwaterloo.ca}
 
\author{Chris Godsil}
\thanks{Support from a National Sciences and Engineering
Council of Canada operating grant is gratefully acknowledged by the second 
author.}
\address{Department of Combinatorics and Optimization,
University of Waterloo, Waterloo, Ontario, Canada N2L 3G1}
\email{cgodsil@uwaterloo.ca}


\begin{abstract}
In 1989, Vaughan Jones introduced spin models and showed that they
could be used to form link invariants in two different ways---by
constructing representations of the braid group, or by constructing
partition functions.  These spin models were subsequently generalized
to so-called 4-weight spin models by Bannai and Bannai; these could be
used to construct partition functions, but did not lead to braid group
representations in any obvious way.  Jaeger showed that spin models
were intimately related to certain association schemes.  
Yamada gave a construction of a symmetric spin model on $4n$ vertices
from each 4-weight spin model on $n$ vertices.

In this paper we build on recent work with Munemasa to
give a different proof to Yamada's result,
and we analyse the structure of the
association scheme attached to this spin model.
\end{abstract}

\maketitle
\section{Introduction}
\label{section_Introduction}

Spin models are a special class of matrices introduced by Jones in
\cite{VJ_Knot} as a tool for creating link invariants.  There are two
strands to their subsequent development that are of interest to us.
First, Jaeger and Nomura showed that all spin models could be realized as
matrices in association schemes (see \cite{N_Alg}).  
Hence spin models have a
combinatorial aspect and, perhaps more importantly, the search for new
spin models was reduced to the search for certain special classes of
association schemes.  (This means that the search space is discrete
rather than continuous.)

The second strand was the development of more general classes of
models, culminating in the four-weight spin models of Bannai and
Bannai \cite{BB_4wt}.  These models are formed from a pair of
matrices; they still provided link invariants, but apparently lacked
the intimate connection to association schemes.

In \cite{CGM}, Munemasa and the present authors developed a new approach
to spin models, based on what we called \textsl{Jones pairs}.  We
showed that these included the four-weight spin models as a special
case.  As a result we were able to show that each four-weight spin
model determines a pair of association schemes.  

In \cite{Yamada2}, Yamada showed that each four-weight spin
model of order $n$ embeds in a very natural way in a spin model of
order $4n$.
We give a complete and different proof to Yamada's result.
In addition, the tools we develop in
Sections~\ref{section_InvJP} to \ref{section_Nom_InvJP} allow us to
analyze the structure of $\Nom{V}$, which was not investigated in
\cite{Yamada2}.

\section{Invertible Jones Pairs}
\label{section_InvJP}

Given two matrices $A$ and $B$ of the same order, 
we use $A \schur B$ to denote their Schur product,
which has
\[
(A \schur B)_{i,j} = A_{i,j} B_{i,j}.
\]
If all entries of $A$ are non-zero, then we say $A$ is
\textsl{Schur invertible} and 
define its \textsl{Schur-inverse}, $\Sinv{A}$, by
\[
\Sinv{A}_{i,j} = \frac{1}{A_{i,j}}.
\]
Equivalently, we have $\Sinv{A} \schur A = J$, 
where $J$ is the matrix of all ones. 

For any $\nbyn$ matrix $C$, we define two linear operators
$X_C$ and $\Delta_C$ as follows:
\[
X_C(M):= CM, \quad \Delta_C(M):=C\schur M \quad
\text{ for all } M \in \mat{n}{\com}.
\]
Given a linear operator $Y$ on $\mat{n}{\com}$,
we use $Y^T$ to denote its adjoint relative to the non-degenerate
bilinear form $tr(M^TN)$ on $\mat{n}{\com}$, and call it the
\textsl{transpose} of $Y$.  It is easy to see that
\begin{equation*}
X^T_C = X_{C^T}, \quad \Delta_C^T=\Delta_C. 
\end{equation*}
A \textsl{Jones pair} is a pair of $\nbyn$ complex matrices $(A,B)$ such
that $X_A$ and $\Delta_B$ are invertible and
\begin{eqnarray}
\label{eqn_JP_1}
\XDX{A}{B}{A} &=& \DXD{B}{A}{B},\\
\label{eqn_JP_2}
\XDX{A}{B^T}{A} &=& \DXD{B^T}{A}{B^T}.
\end{eqnarray}
Note that $X_A$ and $\Delta_B$ are invertible only if
$A$ is invertible and $B$ is Schur
invertible.  It is also easy to observe that $(A,B)$ is a Jones pair
if and only if $(A,B^T)$ is a Jones pair.  Jones pairs are designed to
give representation of braid groups using Jones' construction.  Please
see Section~2 of \cite{CGM} for a description of the construction.

An $\nbyn$ matrix $W$ is a \textsl{type-II matrix} if
\[
W \SinvT{W} = nI.
\]
Note that a type-II matrix is invertible with respect to both matrix
multiplication and the Schur product.  We say that a Jones pair
$(A,B)$ is \textsl{invertible} if $A$ is Schur invertible and $B$ is
invertible.  Theorems~7.1 and 7.2 of \cite{CGM} imply that
a Jones pair $(A,B)$ is invertible if and only if $A$ and $B$
are type-II matrices.

Let $W_1$, $W_2$, $ W_3$ and $ W_4$ be $\nbyn$ complex matrices 
and let $d$ be such that $d^2=n$.
A \textsl{four-weight spin model} is a 5-tuple $(W_1,W_2,W_3,W_4;d)$
that satisfies
\begin{eqnarray}
\label{eqn_4wt_1}
W_3=\SinvT{W_1}, &\quad & W_2=\SinvT{W_4},\\
\label{eqn_4wt_2}
W_1 W_3 = nI,  &\quad & W_2 W_4=nI,\\
\label{eqn_4wt_3}
\sum_{h=1}^{n} (W_1)_{k,h}(W_1)_{h,i}(W_4)_{h,j} &=&
d (W_4)_{i,j}(W_1)_{k,i}(W_4)_{k,j},\\
\label{eqn_4wt_4}
\sum_{h=1}^{n} (W_1)_{h,k}(W_1)_{i,h}(W_4)_{j,h} &=&
d (W_4)_{j,i}(W_1)_{i,k}(W_4)_{j,k}.
\end{eqnarray}
From (\ref{eqn_4wt_1}) and (\ref{eqn_4wt_2}),
we see that both $W_1$ and $W_4$ are type-II matrices and
they determine $W_3$ and $W_2$, respectively.
Furthermore, it is straightforward to verify that 
Equations~(\ref{eqn_4wt_3}) and (\ref{eqn_4wt_4}) 
are equivalent to Equations~(\ref{eqn_JP_1}) and (\ref{eqn_JP_2})
when $W_1=dA$ and $W_4=B$.

Jaeger showed in \cite{J_4wt} that $(A,B)$ and $(C,B)$ are invertible
Jones pairs if and only if $C=DA\inv{D}$ for some invertible diagonal
matrix $D$.  We say that these two invertible Jones pairs are
\textsl{odd-gauge equivalent}.  Proposition~7 of \cite{J_4wt} states that for
every invertible Jones pair $(A,B)$, there exists an invertible
diagonal matrix $D$ such that $DA\inv{D}$ is symmetric.  Since
odd-gauge equivalent invertible Jones pairs give the same link
invariants, we suffer no loss by considering only invertible Jones
pairs whose first matrix is symmetric.

\section{Nomura Algebras}
\label{section_Nom}

We start this section by defining the \textsl{Nomura algebras} $\N{A}{B}$ and
$\dualN{A}{B}$ of a pair of $\nbyn$ matrices.  When $A$ is a type-II
matrix and $B=\Sinv{A}$, our construction gives the Nomura algebras
discussed in \cite{JMN} and \cite{N_Alg}.  The definitions here are
taken from \cite{CGM}.

Let $A$ and $B$ be $\nbyn$ matrices, let $e_1, \ldots, e_n$ be the
standard basis vectors in $\com^n$ and form the $n^2$ column vectors
\[
\Evector{A}{i}{B}{j} \quad \text{for } \range{i,j}{1}{n}.
\]
We define $\N{A}{B}$ to be the set of matrices of which
$\Evector{A}{i}{B}{j}$ is an eigenvector, for all $\range{i,j}{1}{n}$.
This set of matrices is closed under matrix multiplication and
contains the identity matrix $I_n$.

For each matrix $M \in \N{A}{B}$, we use $\T{A}{B}(M)$ to denote
the $\nbyn$ matrix that satisfies
\[
M (\Evector{A}{i}{B}{j}) = \T{A}{B}(M)_{i,j} (\Evector{A}{i}{B}{j}).
\]
We view $\T{A}{B}$ as a linear map from $\N{A}{B}$ to 
$\mat{n}{\com}$ and 
we use $\dualN{A}{B}$ to denote the image of $\N{A}{B}$.
By the definition of $\T{A}{B}$, we have
\[
\T{A}{B}(MN) = \T{A}{B}(M) \schur \T{A}{B}(N).
\]
Consequently the space $\dualN{A}{B}$ is closed under the Schur product.
Since $I_n \in \N{A}{B}$, the matrix $\T{A}{B}(I_n)=J_n$ belongs to
$\dualN{A}{B}$. 
We conclude that $\dualN{A}{B}$ is a commutative algebra
with respect to the Schur product.

If $A$ is invertible, then the columns of $A$ are linearly independent.
Further if $B$ is Schur invertible, 
then for any $j$
\[
\{\Evector{A}{1}{B}{j}, \Evector{A}{2}{B}{j}, \ldots,
\Evector{A}{n}{B}{j}\}
\]
is a basis of $\com^n$.
In this case, the map $\T{A}{B}$ is an isomorphism from 
$\N{A}{B}$, as an algebra with respect to the matrix multiplication,
to $\dualN{A}{B}$, as an algebra with respect to the Schur product. 
We conclude from the commutativity of $\dualN{A}{B}$ 
that $\N{A}{B}$ is commutative with respect to matrix multiplication.

The following result is called the \textsl{Exchange Lemma}.
It will serve as a powerful tool in
Sections~\ref{section_B} and \ref{section_V}.  
The proof of Theorem~\ref{thm_T_AB_XDX} also demonstrates the usefulness
of this lemma.

\begin{lemma}[\cite{CGM}, Lemma~5.1][Exchange]
\label{lem_Exchange}
If $A,B,C,Q,R,S \in \mat{n}{\com}$ then
\[
\XDX{A}{B}{C}=\DXD{Q}{R}{S}
\]
if and only if 
\[
\XDX{A}{C}{B}=\DXD{R}{Q}{S^T}.\qed
\]
\end{lemma}

\begin{theorem}
\label{thm_T_AB_XDX}
If $A$ and $B$ are $\nbyn$ type-II matrices, then
the following are equivalent:
\begin{enumerate}
\item
\label{eqn_T_AB_XDX_a}
$R\in \N{A}{B}$ and $S=\T{A}{B}(R)$.
\item
\label{eqn_T_AB_XDX_b}
$\XDX{R}{B}{A}=\DXD{B}{A}{S}.$
\item
\label{eqn_T_AB_XDX_c}
$\XDX{R}{A}{B}=\DXD{A}{B}{S^T}.$
\item
\label{eqn_T_AB_XDX_d}
$\DXD{B^T}{\SinvT{B}}{nR}=\XDX{S^T}{\SinvT{A}}{A^T}$.
\item
\label{eqn_T_AB_XDX_e}
$\DXD{\SinvT{A}}{A^T}{nR^T}=\XDX{S}{B^T}{\SinvT{B}}$.
\end{enumerate}
\end{theorem}
{\sl Proof.}
The equivalence of (\ref{eqn_T_AB_XDX_a}) and (\ref{eqn_T_AB_XDX_b})
follows from Theorem~6.2 of \cite{CGM}.
 
Applying the Exchange Lemma to (\ref{eqn_T_AB_XDX_b})
gives (\ref{eqn_T_AB_XDX_c}), which is equivalent to
\begin{equation}\label{eqn_Pf_T_AB_XDX_0}
\DXD{\Sinv{A}}{R}{A} = \XDX{B}{S^T}{\inv{B}}.
\end{equation}
Applying the Exchange Lemma to Equation~(\ref{eqn_Pf_T_AB_XDX_0}) again, 
we get
\[
\DXD{R}{\Sinv{A}}{A^T} = \XDX{B}{\inv{B}}{S^T}.
\]
Now we have $\inv{B}=n^{-1} \SinvT{B}$ and $\Sinv{A} = n \invT{A}$
because $A$ and $B$ are type-II matrices.  The above equation becomes
\[
\DXD{R}{n\invT{A}}{A^T} = \XDX{B}{n^{-1}\SinvT{B}}{S^T},
\]
which leads to
\begin{equation}\label{eqn_Pf_T_AB_XDX_1}
\DXD{B^T}{\inv{B}}{nR} = \XDX{S^T}{\SinvT{A}}{n^{-1} A^T}.
\end{equation}
We get (\ref{eqn_T_AB_XDX_d}) after multiplying both sides of 
Equation~(\ref{eqn_Pf_T_AB_XDX_1})
by $n$ and replacing $n\inv{B}$ by $\SinvT{B}$.

Taking the transpose of both sides of Equation~(\ref{eqn_Pf_T_AB_XDX_1})
gives
\[
\DXD{nR}{\invT{B}}{B^T} = \XDX{n^{-1}A}{\SinvT{A}}{S}
\]
and
\[
\DXD{A^T}{\SinvT{A}}{nR} = \XDX{S}{\SinvT{B}}{B^T}.
\]
We get (\ref{eqn_T_AB_XDX_e}) after applying
the Exchange Lemma to the above equation.\qed

Now we state an easy consequence of
Theorem~\ref{thm_T_AB_XDX}~(\ref{eqn_T_AB_XDX_b}). 
\begin{corollary}[\cite{CGM}, Lemma~10.2]
\label{cor_Sinv_AB}
Let $A$ and $B$ be $\nbyn$ type-II matrices.
If $R \in \N{A}{B}$ then
\[
R^T \in \N{\Sinv{A}}{\Sinv{B}}, \quad \text{and } \quad
\T{\Sinv{A}}{\Sinv{B}}(R^T) = \T{A}{B}(R).\qed
\]
\end{corollary}

\section{Nomura Algebras of a Type-II Matrix}
\label{section_Nom_II}

When $A$ is a type-II matrix and $B=\Sinv{A}$, existing papers such as
\cite{JMN} use $\Nom{A}$, $\dualNom{A}$ and $\Theta_A$ to denote
$\N{A}{B}$ , $\dualN{A}{B}$ and $\T{A}{B}$, respectively.  The algebra
$\Nom{A}$ is called the Nomura algebra of $A$.  We now present some
results on $\Nom{A}$ due to Jaeger, Matsumoto and Nomura \cite{JMN}
which we will use later.

When $B=\Sinv{A}$, Condition~\ref{thm_T_AB_XDX}~(\ref{eqn_T_AB_XDX_e})
becomes 
\[
\DXD{\SinvT{A}}{A^T}{nR^T} = \XDX{S}{\SinvT{A}}{A^T},
\]
and it implies
\begin{equation}
\label{eqn_Duality_A}
\Theta_{A^T}(S)=\Theta_{A^T}(\Theta_A(R)) = n R^T.
\end{equation}
We conclude that if $R \in \Nom{A}$ then 
$\Theta_A(R) \in \Nom{A^T}$ and $R^T \in \dualNom{A^T}$.
Hence 
\[
\dualNom{A} \subseteq \Nom{A^T} \quad
\text{and}\quad
\dim(\Nom{A}) =\dim(\Nom{A}^T)\leq \dim(\dualNom{A^T}).
\]
Similarly $A^T$ is also a type-II matrix, so 
\[
\dualNom{A^T} \subseteq \Nom{A} \quad
\text{and}\quad
\dim(\Nom{A^T}) \leq \dim(\dualNom{A}).
\]
Therefore $\dualNom{A}=\Nom{A^T}$ and $\dualNom{A^T}=\Nom{A}$,
which implies that $\Nom{A}$ and $\Nom{A^T}$ are closed
under both matrix multiplication and the Schur product.
It also implies that $\Nom{A}=\dualNom{A^T}$ is closed under the transpose.
Since $A$ is invertible and $\Sinv{A}$ is Schur invertible,
the map $\Theta_A$ is an isomorphism from $\Nom{A}$ to $\dualNom{A}$.
Hence $\Nom{A}$ is commutative with respect to matrix multiplication.
In summary, the algebra $\Nom{A}$ is commutative with respect to matrix
multiplication, is also closed under the transpose and the Schur product,
and contains $I$ and $J$.  In other words, $\Nom{A}$ is a 
\textsl{Bose-Mesner algebra}.

We now investigate the properties of the map $\Theta_A$.
Let $M$ and $N$ be matrices in $\Nom{A}$.
Since $\Theta_{A^T}:\Nom{A^T} \rightarrow \Nom{A}$ is an isomorphism,
there exist $M'$ and $N'$ in $\Nom{A^T}$ such that
$\Theta_{A^T}(M')=M$ and  $\Theta_{A^T}(N') = N$.
Hence
\begin{eqnarray*}
\Theta_A(M\schur N) &=& \Theta_A(\Theta_{A^T}(M') \schur \Theta_{A^T}(N'))\\
&=& \Theta_A(\Theta_{A^T}(M'N'))
\end{eqnarray*}
which equals $n (M'N')^T$ by Equation~(\ref{eqn_Duality_A}).
Since 
\[
\Theta_A(M)=\Theta_A(\Theta_{A^T}(M'))= nM'^T
\] 
and $\Theta_A(N)=nN'^T$, we have
\begin{eqnarray*}
\Theta_A(M \schur N) 
&=& \frac{1}{n} (n N'^T)(nM'^T)\\
&=& \frac{1}{n} \Theta_A(N) \Theta_A(M)\\
&=& \frac{1}{n} \Theta_A(M) \Theta_A(N),
\end{eqnarray*}
the last equality results from the commutativity of $\dualNom{A}$.
Now we conclude that $\Theta_A$ swaps matrix multiplication
with the Schur product. 

Furthermore, applying $\frac{1}{n}\Theta_A$ to the two rightmost terms of 
Equation~(\ref{eqn_Duality_A}) gives
\[
\frac{1}{n} \Theta_A\left( \Theta_{A^T}(\Theta_A(R)) \right)=
\Theta_A(R^T).
\]
It follows from Equation~(\ref{eqn_Duality_A}) that the left-hand side
equals $\Theta_A(R)^T$.
Thus $\Theta_A$ and the transpose commute.
From Corollary~\ref{cor_Sinv_AB}, we see that 
\[
\Theta_{\Sinv{A}}(R)= \Theta_A(R)^T.
\]
Also note that by Equation~(\ref{eqn_Duality_A}), we have 
\[
\Theta_A(J) = \Theta_A(\Theta_{A^T}(I))=nI.
\]
We call $\Theta_A$ a \textsl{duality map} from $\Nom{A}$ to $\Nom{A^T}$ and say
that these two Bose-Mesner algebras form a \textsl{formally dual pair}.  If
$\Nom{A}=\Nom{A^T}$ and
$\Theta_A = \Theta_{A^T}$, we say that it is \textsl{formally self-dual}.

A \textsl{spin model} is an $\nbyn$ matrix $W$ such that
$(W,W,\Sinv{W},\Sinv{W};d)$ is a four-weight spin model, for $d^2=n$.
It follows from Section~9 of \cite{CGM} that $W$ is a spin model if
and only if $(d^{-1} W, \Sinv{W})$ is an invertible Jones pair.  In
\cite{JMN}, Jaeger, Matsumoto and Nomura gave the following
characterization of a spin model $W$ using its Nomura algebra
$\Nom{W}$.

\begin{theorem}[\cite{JMN}, Theorem~11]
\label{thm_Spin}
Suppose $W$ is a type-II matrix.  Then $W \in \Nom{W}$ if and only if
$cW$ is a spin model for some non-zero scalar $c$.  In this case,
\[
\Nom{W} = \Nom{W^T}
\]
is a formally self-dual Bose-Mesner algebra with duality map
$\Theta_W=\Theta_{W^T}$.\qed
\end{theorem}

\section{Nomura Algebras of an Invertible Jones Pair}
\label{section_Nom_InvJP}

We study the relation among the different Nomura algebras of 
an invertible Jones pair.
\begin{theorem}[\cite{Ban1}, Theorem~3]
\label{thm_InvJP_Nom}
If $(A,B)$ is an invertible Jones pair, then
\[
\Nom{A} = \Nom{A^T} = \Nom{B} = \Nom{B^T},
\]
the duality maps satisfy
$\Theta_A = \Theta_{A^T}$ and $\Theta_B = \Theta_{B^T}$.\qed
\end{theorem}
Bannai, Guo and Huang \cite{Ban1} proved this result for four-weight
spin models, which are equivalent to invertible Jones pairs.
For an alternate proof using the Nomura algebras of $A$ and $B$, see
Section~10 of \cite{CGM}.

Let $A$ and $B$ be type-II matrices.
We see from Theorem~\ref{thm_T_AB_XDX}~(\ref{eqn_T_AB_XDX_a})
and (\ref{eqn_T_AB_XDX_b}) that $(A,B)$ is an invertible
Jones pair if and only if 
$A \in \N{A}{B} \cap \N{A}{B^T}$, $\T{A}{B}(A)=B$ and $\T{A}{B^T}(A)=B^T$.
The next two results provide some insights to the relations among
$\N{A}{B}$, $\dualN{A}{B}$ and $\Nom{A}$. 
\begin{theorem}[\cite{CGM}, Theorem~10.3]
\label{thm_T_AB}
Let $A$ and $B$ be $\nbyn$ type-II matrices.
If $F\in \Nom{A}$, $G \in \N{A}{B}$ and $H\in \Nom{B}$, then
$F \schur G$, and $G \schur H$ belong to $\N{A}{B}$ and 
\begin{eqnarray*}
\T{A}{B}(F \schur G) &=& n^{-1} \Theta_A(F)\ \T{A}{B}(G),\\
\T{A}{B}(G \schur H) &=& n^{-1} \T{A}{B}(G)\ \Theta_B(H)^T. 
\end{eqnarray*}
\qed
\end{theorem}

\begin{theorem}[\cite{CGM}, Theorem~10.4]
\label{thm_TA_TB}
Let $A$ and $B$ be $\nbyn$ type-II matrices.
If $F,G \in \N{A}{B}$, then $F \schur G^T \in \Nom{A} \cap \Nom{B}$ and 
\begin{eqnarray}
\nonumber
\Theta_A(F \schur G^T) &=& n^{-1} \T{A}{B}(F)\ \T{A}{B}(G)^T,\\
\Theta_B(F \schur G^T) &=& n^{-1} \T{A}{B}(F)^T\ \T{A}{B}(G). 
\label{eqn_thm_TA_TB}
\end{eqnarray}
\qed
\end{theorem}

We list two consequences of Theorems~\ref{thm_T_AB} and \ref{thm_TA_TB}.
\begin{theorem}[\cite{CGM}, Theorem~10.6]
\label{thm_NA_NAB_GH}
Let $A$ and $B$ be $\nbyn$ type-II matrices.  
If $\N{A}{B}$ contains a Schur invertible matrix $G$ and $H=\T{A}{B}(G)$,
then
\[
\N{A}{B}=G \schur \Nom{A},
\qquad
\dualN{A}{B}H^T=\Nom{A^T}.\qed
\]
\end{theorem}

\begin{corollary}[\cite{CGM}, Corollary~10.9]
\label{cor_TA_TB}
If $(A,B)$ is an invertible Jones pair, then
\[
\Theta_B(M)^T = \inv{B}\Theta_A(M)B
\]
for all $M \in \Nom{A}$.\qed
\end{corollary}

Now we present an important application of Theorems~\ref{thm_T_AB} and
\ref{thm_TA_TB}, which implies that the Nomura algebras $\Nom{A}$,
$\N{A}{B}$ and $\dualN{A}{B}$ have the same dimension.

\begin{theorem}
\label{thm_NA_NAB}
Let $(A,B)$ be an invertible Jones pair.  Then
\[
\N{A}{B} = A \schur \Nom{A}, \qquad
\dualN{A}{B}B^T = \Nom{A}
\]
and
\[
\dualN{A}{B} = (\dualN{A}{B^T})^T.
\]
\end{theorem}
{\sl Proof.}
We get the first equality by letting $G=A$ in Theorem~\ref{thm_NA_NAB_GH}.
Since $B=\T{A}{B}(A)$, we have
\[
\dualN{A}{B}B^T=\Nom{A^T}.
\]
By Theorem~\ref{thm_InvJP_Nom}, we have $\Nom{A^T}=\Nom{A}$ and
hence the second equality holds.

If we replace $B$ by $B^T$ in the above equality, then
we get 
\[
\dualN{A}{B^T}B=\Nom{A^T}.
\]
Since multiplication by $B$ is injective, the dimensions of $\Nom{A}=\Nom{A^T}$
and $\dualN{A}{B^T}$ are equal.  Now we let $G$ equal $A$ and
replace $B$ by $B^T$ in Equation~(\ref{eqn_thm_TA_TB}).
We get
\[
\dualNom{B^T} \subseteq (\dualN{A}{B^T})^T B^T.
\]
By Theorem~\ref{thm_InvJP_Nom}, $\Nom{A}=\Nom{B}=\dualNom{B^T}$.
Since $\Nom{A}$ and $\dualN{A}{B^T}$ have the same dimension,
we have 
\[
\Nom{A}=(\dualN{A}{B^T})^T B^T.
\]
Thus $\dualN{A}{B} B^T=(\dualN{A}{B^T})^T B^T$, which leads to
the last equality of the theorem.\qed

\begin{corollary}
\label{cor_NAB}
Let $(A,B)$ be an invertible Jones pair.
Then 
\[
\N{A}{B} = \N{A}{B^T}.
\]
Moreover, if $A$ is symmetric, then
\[
\N{A}{B} = (\N{A}{B})^T.
\]
\end{corollary}
{\sl Proof.}
Applying Theorem~\ref{thm_NA_NAB} to the invertible Jones pairs 
$(A,B)$ and $(A,B^T)$ gives
\[
\N{A}{B} = A \schur \Nom{A} = \N{A}{B^T}.
\]
Using the same equation, we have
$\N{A}{B}^T= A^T \schur \Nom{A}^T$.
Since $\Nom{A}$ is closed under the transpose and $A$ is symmetric, we
conclude that $\N{A}{B}=\N{A}{B}^T$.\qed

\section{A Bose-Mesner Algebra of order $4n$}
\label{section_B}

From now on, we assume that $(A,B)$ is an invertible Jones pair and
$A$ is symmetric.
\begin{lemma}
\label{lem_cond1}
For each $H$ in $\N{A}{B}$, there exists a unique matrix $K$ in
$(\N{A}{B^T})^T$ such that
\begin{equation}
\label{cond1}
\T{A}{B}(H)=\T{A}{B^T}(K^T)^T.
\end{equation}
\end{lemma}
{\sl Proof.}
Existence follows directly from the last equality in
Theorem~\ref{thm_NA_NAB}, while uniqueness holds because
$\Theta_{A,B^T}$ is an isomorphism.\qed 
Given any matrix $H$ in $\N{A}{B}$, 
we say that the unique $K$ in $\dualN{A}{B}$
satisfying Equation~(\ref{cond1}) is {\sl paired} with $H$.
\begin{lemma}
\label{lem_cond2}
For each $H$ in $\N{A}{B}$, $K$ in $\dualN{A}{B}$ is paired with $H$
if and only if $K^T$ is paired with $H^T$.
Moreover we have
\begin{equation}
\label{cond2}
\Theta_A(H \schur A)=\Theta_{B^T}(K^T \schur A).
\end{equation}
\end{lemma}
{\sl Proof.}
Multiplying each side of Equation~(\ref{cond1}) by
$n^{-1}\T{A}{B}(A)^T=n^{-1}\T{A}{B^T}(A)$ gives
\[
n^{-1}\T{A}{B}(H)\T{A}{B}(A)^T=n^{-1}\T{A}{B^T}(K^T)^T\T{A}{B^T}(A).
\]
We apply Theorem~\ref{thm_TA_TB} to both sides of the
above equation to get
\[
\Theta_A(H \schur A^T)=\Theta_{B^T}(K^T \schur A^T).
\]
Since $A$ is symmetric,
we see that Equation~(\ref{cond1}) is equivalent to
Equation~(\ref{cond2}).

In addition, taking the transpose of both sides gives
\[
\Theta_A(H^T \schur A)=\Theta_{B^T}(K \schur A).
\]
Therefore $H$ and $K$ satisfy Equation~(\ref{cond1}) if and only if
$H^T$ and $K^T$ satisfy Equation~(\ref{cond1}).
\qed

For any $F \in \Nom{A}$ and $H,G \in \N{A}{B}$, we define the $4n
\times 4n$ matrix $\M(F,G,H)$ to be
\[
\begin{pmatrix}
\Theta_A(F)+H & \Theta_A(F)-H & \T{A}{B}(G) & \T{A}{B}(G)\\
\Theta_A(F)-H & \Theta_A(F)+H & \T{A}{B}(G) & \T{A}{B}(G)\\
\T{A}{B}(G^T)^T & \T{A}{B}(G^T)^T & 
\Theta_{\Sinv{B}}(F)+K & \Theta_{\Sinv{B}}(F)-K\\
\T{A}{B}(G^T)^T & \T{A}{B}(G^T)^T & 
\Theta_{\Sinv{B}}(F)-K & \Theta_{\Sinv{B}}(F)+K\\
\end{pmatrix},
\]
where $K$ is paired with $H$.
We consider the space
\begin{equation}
\label{eqn_B}
\BM:=\{\M(F,G,H) : F \in \Nom{A} \text{ and } H, G \in \N{A}{B}\}.
\end{equation}
Now we show that $\BM$ is a Bose-Mesner algebra.
It turns out that $\BM$ contains the $4n \times 4n$ type-II matrix $V$ 
defined at the beginning of Section~\ref{section_V} and 
it is a subscheme of $\Nom{V}$.
This leads to the main result of this paper which says that
$V$ is a spin model if and only if $(A,B)$ is an invertible Jones pair.

To convince ourselves that $\BM$ is a Bose-Mesner algebra, we need to
check that $\BM$ contains the identity matrix $I_{4n}$ and 
the matrix of all ones $J_{4n}$; 
it is closed under the transpose;
it is a commutative algebra with respect to matrix multiplication;
it is closed under the Schur product.

\begin{lemma}
\label{lem_calB_IJ}
The vector space $\BM$ contains $I_{4n}$ and $J_{4n}$.
\end{lemma}
{\sl Proof.}
The matrix $K$ that is paired with $\frac{1}{2}I_n$ satisfies
\[
\T{A}{B^T}(K^T)^T = \T{A}{B}(\frac{1}{2}I_n)=\frac{1}{2}J_n.
\]
Since $\T{A}{B^T}$ is an isomorphism, we conclude that $K=\frac{1}{2}I_n$.
Note that $\Theta_A(\frac{1}{2n}J_n) = \frac{1}{2}I_n$.
Thus $\M(\frac{1}{2n}J_n,\0,\frac{1}{2}I_n)=I_{4n}$ belongs to $\BM$.

Since $\Theta_A(I_n)=\T{A}{B}(I_n)=J_n$, the matrix
$\M(I_n,I_n,\0)=J_{4n}$ belongs to $\BM$.\qed

\begin{lemma}
\label{lem_calB_transpose}
The vector space $\BM$ is closed under transpose.
\end{lemma}
{\sl Proof.}
Let $\M(F,G,H) \in \BM$. 
Now $\M(F,G,H)^T$ equals
\begin{eqnarray*}
\begin{pmatrix}
\Theta_A(F)^T+H^T & \Theta_A(F)^T-H^T & \T{A}{B}(G^T) & \T{A}{B}(G^T)\\
\Theta_A(F)^T-H^T & \Theta_A(F)^T+H^T & \T{A}{B}(G^T) & \T{A}{B}(G^T)\\
\T{A}{B}(G)^T & \T{A}{B}(G)^T & 
\Theta_{\Sinv{B}}(F)^T+K^T & \Theta_{\Sinv{B}}(F)^T-K^T\\
\T{A}{B}(G)^T & \T{A}{B}(G)^T & 
\Theta_{\Sinv{B}}(F)^T-K^T & \Theta_{\Sinv{B}}(F)^T+K^T\\
\end{pmatrix}.
\end{eqnarray*}
Since $\N{A}{B}$ is closed under the transpose, the matrices
$G^T$ and $H^T$ belong to $\N{A}{B}$.
It follows from Lemma~\ref{lem_cond2} that $K^T$ is paired with $H^T$. 
Moreover, $\Theta_A(F)^T = \Theta_A(F^T)$.
As a result we conclude that 
\[
\M(F,G,H)^T = \M(F^T,G^T,H^T),
\]
and the vector space $\BM$ is closed under the transpose.\qed

\begin{lemma}
\label{lem_calB_mult}
The vector space $\BM$ is a commutative algebra under matrix multiplication.
\end{lemma}
{\sl Proof.}
Let $M=\M(F,G,H)$ and $M_1=\M(F_1,G_1,H_1)$ be any matrices in $\BM$.

By Theorem~\ref{thm_TA_TB}, we have 
\[
\T{A}{B}(G)\T{A}{B}(G_1^T)^T = n \Theta_A(G \schur G_1).
\]
Hence the top left $2n \times 2n$ block of $M M_1$ equals
\[
\begin{pmatrix}
2n\Theta_A(F \schur F_1 + G \schur G_1)+2HH_1 &
2n\Theta_A(F \schur F_1 + G \schur G_1)-2HH_1 \\
2n\Theta_A(F \schur F_1 + G \schur G_1)-2HH_1 &
2n\Theta_A(F \schur F_1 + G \schur G_1)+2HH_1 \\
\end{pmatrix}.
\]
Similarly, by Theorem~\ref{thm_TA_TB}
\begin{eqnarray*}
\T{A}{B}(G^T)^T\T{A}{B}(G_1) 
&=& n \Theta_{B}(G^T \schur G_1^T)\\
&=& n \Theta_{B}(G \schur G_1)^T\\
&=& n \Theta_{\Sinv{B}}(G \schur G_1).
\end{eqnarray*}
Consequently the bottom right $2n \times 2n$ block of $M M_1$ equals
\[
\begin{pmatrix}
2n \Theta_{\Sinv{B}}(F \schur F_1 + G \schur G_1) + 2K K_1 &
2n \Theta_{\Sinv{B}}(F \schur F_1 + G \schur G_1) - 2K K_1 \\
2n \Theta_{\Sinv{B}}(F \schur F_1 + G \schur G_1) - 2K K_1 &
2n \Theta_{\Sinv{B}}(F \schur F_1 + G \schur G_1) + 2K K_1 \\
\end{pmatrix},
\]
where $K$ and $K_1$ are paired with $H$ and $H_1$, respectively.
Now we need to show that $K K_1$ is paired with $H H_1$.
From Equation~(\ref{cond1}), we have
\[
\T{A}{B}(H)=\T{A}{B^T}(K^T)^T \quad \text{and} \quad
\T{A}{B}(H_1)=\T{A}{B^T}(K_1^T)^T.
\]
Therefore
\begin{eqnarray*}
\T{A}{B}(H H_1) &=& 
\T{A}{B}(H) \schur \T{A}{B}(H_1) \\
&=& \T{A}{B^T}(K^T)^T\schur \T{A}{B^T}(K_1^T)^T\\
&=& \T{A}{B^T}(K^T K_1^T)^T.
\end{eqnarray*}
Since $\N{A}{B^T}$ is commutative with respect to matrix multiplication,
\[
\T{A}{B}(H H_1) = \T{A}{B^T}((KK_1)^T)^T.
\]

We now consider the top right $2n \times 2n$ block of $M M_1$.
Note that
\begin{eqnarray*}
&&2\Theta_A(F)\T{A}{B}(G_1) + 2\T{A}{B}(G)\Theta_{\Sinv{B}}(F_1)\\
&=&2\Theta_A(F)\T{A}{B}(G_1) + 2\T{A}{B}(G)\Theta_{B}(F_1)^T.
\end{eqnarray*}
Applying Theorem~\ref{thm_T_AB} to each term, we get
\[
2n \T{A}{B}(F \schur G_1 + G \schur F_1).
\]
Thus the top right $2n \times 2n$ block of $M M_1$ is
\[
\begin{pmatrix}
2n \T{A}{B}(F \schur G_1 + G \schur F_1) &
2n \T{A}{B}(F \schur G_1 + G \schur F_1) \\
2n \T{A}{B}(F \schur G_1 + G \schur F_1) &
2n \T{A}{B}(F \schur G_1 + G \schur F_1) \\
\end{pmatrix}.
\]
Consider the bottom left $2n \times 2n$ block of $M M_1$,
we have
\begin{eqnarray*}
&&2\T{A}{B}(G^T)^T\Theta_A(F_1)+2\Theta_{\Sinv{B}}(F)\T{A}{B}(G_1^T)^T\\
&=&2\T{A}{B}(G^T)^T\Theta_A(F_1)+2\Theta_{B}(F)^T\T{A}{B}(G_1^T)^T.
\end{eqnarray*}
Since each of $\Theta_A$ and $\Theta_B$ commutes with the transpose,
the above expression becomes
\[
2\T{A}{B}(G^T)^T\Theta_A(F_1^T)^T+2\Theta_{B}(F^T)\T{A}{B}(G_1^T)^T,
\]
which equals
\[
2n \T{A}{B}(F_1^T \schur G^T + G_1^T \schur F^T)^T
\] 
by Theorem~\ref{thm_T_AB}.
Hence the bottom left $2n \times 2n$ block of $M M_1$ is 
\[
\begin{pmatrix}
2n \T{A}{B}(F_1^T \schur G^T + G_1^T \schur F^T)^T &
2n \T{A}{B}(F_1^T \schur G^T + G_1^T \schur F^T)^T \\
2n \T{A}{B}(F_1^T \schur G^T + G_1^T \schur F^T)^T &
2n \T{A}{B}(F_1^T \schur G^T + G_1^T \schur F^T)^T \\
\end{pmatrix}.
\]

Now we conclude that 
\[
M M_1 =  \M(2nF \schur F_1+ 2nG \schur G_1,
2nF\schur G_1+2nG \schur F_1, 2H H_1)
\]
belongs to $\BM$.

It follows from the commutativity of $\N{A}{B}=\N{A}{B^T}$ that 
$H H_1=H_1 H$ and $KK_1=K_1K$.
Therefore all four $2n \times 2n$ blocks of $M M_1$ 
remain unchanged after swapping
$F$ with $F_1$, $G$ with $G_1$, $H$ with $H_1$ and $K$ with $K_1$.
Consequently the matrices $M$ and $M_1$ commute.\qed

\begin{lemma}
\label{lem_calB_schur}
The algebra $\BM$ is closed under the Schur product.
\end{lemma}
{\sl Proof.}
Let $M=\M(F,G,H)$ and $M_1=\M(F_1,G_1,H_1)$ be two matrices in $\BM$.
We want to write $M\schur M_1$ as $\M(F',G',H')$, for some
$F'$ in $\Nom{A}$ and $G'$ and $H'$ in $\N{A}{B}$.
If we divide $M\schur M_1$ into sixteen $\nbyn$ blocks naturally, then
the $(1,1)$- and $(2,2)$-blocks of $M\schur M_1$ are equal to
\begin{eqnarray*}
&&\Theta_A(F)\schur \Theta_A(F_1) + H \schur H_1 
+ \Theta_A(F)\schur H_1 + \Theta_A(F_1) \schur H \\
&=&\left( \Theta_A(F F_1) + H \schur H_1 \right)
+ \left( \Theta_A(F)\schur H_1 + \Theta_A(F_1) \schur H \right).
\end{eqnarray*}
The $(1,2)$- and $(2,1)$-blocks of $M \schur M_1$ are equal to
\[
\left( \Theta_A(F F_1) + H \schur H_1 \right)
- \left( \Theta_A(F)\schur H_1 + \Theta_A(F_1) \schur H \right).
\]
The $(3,3)$- and $(4,4)$-blocks of $M \schur M_1$ are equal to
\[
\left( \Theta_{\Sinv{B}}(F F_1) + K \schur K_1 \right)
+ \left( \Theta_{\Sinv{B}}(F)\schur K_1 + \Theta_{\Sinv{B}}(F_1)\schur
  K\right).
\]
The $(3,4)$- and $(4,3)$-blocks of $M \schur M_1$ are equal to
\[
\left( \Theta_{\Sinv{B}}(F F_1) + K \schur K_1 \right)
- \left( \Theta_{\Sinv{B}}(F)\schur K_1 +\Theta_{\Sinv{B}}(F_1)\schur 
  K \right).
\]
To determine $F'$,
we need to show that there exists $\hat{F} \in \Nom{A}$ such that
\[
H \schur H_1 = \Theta_A(\hat{F})\quad \text{and} \quad
K \schur K_1 = \Theta_{\Sinv{B}}(\hat{F}),
\]
and $F'=F F_1+\hat{F}$.
Now the matrix $K$ is paired with $H$.
Right-multiplying both sides of Equation~(\ref{cond1}) by
$\SinvT{B}$ yields
\[
\T{A}{B}(H) \SinvT{B} = \T{A}{B^T}(K^T)^T \SinvT{B},
\]
which is rewritten as
\[
\T{A}{B}(H) \T{A}{B}(\inv{A})^T =
\T{A}{B^T}(K^T)^T \T{A}{B^T}(\inv{A}).
\]
Since $\inv{A} = \frac{1}{n}\SinvT{A}$, the above
equation is equivalent to
\[
\frac{1}{n}\T{A}{B}(H) \T{A}{B}(\SinvT{A})^T =
\frac{1}{n}\T{A}{B^T}(K^T)^T \T{A}{B^T}(\SinvT{A})
\]
Applying Theorem~\ref{thm_TA_TB} to each side, we get
\[
\Theta_A(H \schur \Sinv{A}) 
= \Theta_{B^T}(K^T \schur \Sinv{A}).
\]
Applying Corollary~\ref{cor_Sinv_AB} to the right-hand side, we get
\begin{eqnarray*}
\Theta_A(H \schur \Sinv{A}) &=& \Theta_{\SinvT{B}}(K \schur \SinvT{A})\\
&=& \Theta_{\SinvT{B}}(K \schur \Sinv{A}).
\end{eqnarray*}
Similarly, $M_1 \in \BM$.  By Lemma~\ref{lem_cond2},
the matrices $H_1$ and $K_1$ satisfy Equation~(\ref{cond2})
\[
\Theta_A(H_1 \schur A)=\Theta_{B^T}(K_1^T \schur A)
=\Theta_{\SinvT{B}}(K_1 \schur A^T).
\]
Since $A$ is symmetric,
\[
\Theta_A(H \schur \Sinv{A})\ \Theta_A(H_1 \schur A) =
\Theta_{\SinvT{B}}(K \schur \Sinv{A})\ \Theta_{\SinvT{B}}(K_1 \schur A)
\]
and 
\[
\Theta_A(H \schur \Sinv{A} \schur H_1 \schur A) =
\Theta_{\SinvT{B}}(K \schur \Sinv{A} \schur K_1 \schur A),
\]
which simplifies to
\[
\Theta_A(H \schur H_1) = \Theta_{\SinvT{B}}(K \schur K_1).
\]
If we let $\hat{F} = \frac{1}{n}\Theta_A(H \schur H_1)^T$, then
\begin{eqnarray*}
\Theta_A(\hat{F}) &=&
\frac{1}{n}\Theta_A(\Theta_A(H\schur H_1)^T)\\
&=& \frac{1}{n} \Theta_A(\Theta_A(H \schur H_1))^T.
\end{eqnarray*}
Since $A$ is symmetric, it follows from Equation~(\ref{eqn_Duality_A})
that $\Theta_A(\hat{F}) = H \schur H_1$ and
\begin{eqnarray*}
\Theta_{\Sinv{B}}(\hat{F}) &=& 
\frac{1}{n}\Theta_{\Sinv{B}}(\Theta_{\SinvT{B}}(K \schur K_1)^T)\\
&=& K \schur K_1.
\end{eqnarray*}
As a result we have $F'=F F_1+\frac{1}{n}\Theta_A(H \schur H_1)^T$.

We see from the $(1,1)$- and $(1,2)$-blocks of $M \schur M_1$ that
$H'$ should be equal to $\Theta_A(F)\schur H_1 + \Theta_A(F_1) \schur H$.  We
now need to verify that 
$\Theta_{\Sinv{B}}(F)\schur K_1+ \Theta_{\Sinv{B}}(F_1)\schur K$
is paired with $H'$.
That is,
\begin{align}
\label{eqn_Verify9}
 \T{A}{B} (\Theta_A(F)\schur H_1 &+ \Theta_A(F_1) \schur H)\\
&= \T{A}{B^T}
\left((\Theta_{\Sinv{B}}(F)\schur K_1+
\Theta_{\Sinv{B}}(F_1)\schur K\right)^T)^T.
\nonumber
\end{align}
Applying Theorem~\ref{thm_T_AB} gives
\begin{align*}
\T{A}{B}(\Theta_A(F) \schur H_1 &+ \Theta_A(F_1)\schur H)\\
&= \frac{1}{n}\Theta_A(\Theta_A(F))\T{A}{B}(H_1)
+ \frac{1}{n}\Theta_A(\Theta_A(F_1))\T{A}{B}(H)\\
&= F^T \T{A}{B}(H_1) + F_1^T\T{A}{B}(H).
\end{align*}
By Equation~(\ref{cond1}), the above expression equals
\[
 F^T \T{A}{B^T}(K_1^T)^T + F_1^T\T{A}{B^T}(K^T)^T
= \left( \T{A}{B^T}(K_1^T)F + \T{A}{B^T}(K^T)F_1 \right)^T.
\]
By Equation~(\ref{eqn_Duality_A}), we see that
$F=\inv n\Theta_{B^T}(\Theta_B(F))^T$ and consequently the above expression
is equal to
\[
\left(\frac{1}{n} \T{A}{B^T}(K_1^T)\Theta_{B^T}(\Theta_B(F))^T +
\frac{1}{n} \T{A}{B^T}(K^T)\Theta_{B^T}(\Theta_B(F_1))^T \right)^T.
\]
Applying Theorem~\ref{thm_T_AB} yields
\begin{align*}
 \T{A}{B^T}(K_1^T\schur \Theta_B(F))^T
&+ \T{A}{B^T}(K^T \schur \Theta_B(F_1))^T\\
&= \T{A}{B^T}(K_1^T\schur \Theta_{\Sinv{B}}(F)^T
+ K^T \schur \Theta_{\Sinv{B}}(F_1)^T)^T\\
&= \T{A}{B^T}\left((K_1 \schur \Theta_{\Sinv{B}}(F)
+ K \schur \Theta_{\Sinv{B}}(F_1))^T\right)^T.
\end{align*}
Hence Equation~(\ref{eqn_Verify9}) is satisfied and
$H'=\Theta_A(F)\schur H_1 + \Theta_A(F_1) \schur H$.

Since 
\[
\T{A}{B}(G) \schur \T{A}{B}(G_1) = \T{A}{B}(G G_1)
\]
and
\begin{eqnarray*}
\T{A}{B}(G_1^T)^T \schur \T{A}{B}(G^T)^T 
&=& \T{A}{B}(G_1^T G^T)^T\\
&=& \T{A}{B}((G G_1)^T)^T,
\end{eqnarray*}
the top right $2n \times 2n$ and the bottom left
$2n \times 2n$ blocks of $M \schur M_1$ are
\[
\begin{pmatrix}
\T{A}{B}(GG_1) & \T{A}{B}(GG_1) \\
\T{A}{B}(GG_1) & \T{A}{B}(GG_1) \\
\end{pmatrix}
\]
and
\[
\begin{pmatrix}
\T{A}{B}((GG_1)^T)^T & \T{A}{B}((GG_1)^T)^T \\
\T{A}{B}((GG_1)^T)^T & \T{A}{B}((GG_1)^T)^T \\
\end{pmatrix},
\]
respectively.
We conclude that $G'=G G_1$ and that
\[
M \schur M_1 = \M(FF_1+\frac{1}{n}\Theta_A(H \schur H_1)^T, GG_1,
\Theta_A(F)\schur H_1 + \Theta_A(F_1) \schur H)
\]
belongs to $\BM$.\qed

\begin{theorem}
\label{thm_calB}
The algebra
$\BM$ is a Bose-Mesner algebra whose dimension is three times 
the dimension of $\Nom{A}$.
\end{theorem}
{\sl Proof.}
It follows from Lemmas~\ref{lem_calB_IJ} to \ref{lem_calB_schur} that
$\BM$ is a Bose-Mesner algebra.
By the definition of the matrices in $\BM$, the algebra $\BM$ is the
direct sum of three vector spaces.  The first one consists of matrices
$\M(F,\0,\0)$ for all $F \in \Nom{A}$.  This space is isomorphic to
$\Nom{A}$.  The second vector space consists of matrices
$\M(\0,G,\0)$ for all $G \in \N{A}{B}$.  The third one consists of
matrices $\M(\0,\0,H)$ for all $H \in \N{A}{B}$.  Both the second and
the third vector spaces are isomorphic to $\N{A}{B}$.  By
Theorem~\ref{thm_NA_NAB_GH}, $\Nom{A}$ and $\N{A}{B}$ have the same
dimension.  Therefore the dimension of $\BM$ is three times the
dimension of $\Nom{A}$.\qed

\section{A $4n \times 4n$ Symmetric Spin Model}
\label{section_V}

Let $A$ and $B$ be $\nbyn$ type-II matrices, and assume $A$ is symmetric.
Let $d$ be such that $d^2=n$.
In \cite{Yamada2}, Yamada defined
a symmetric $4n \times 4n$ matrix
\[
V:=
\begin{pmatrix}
dA & -dA & \Sinv{B} & \Sinv{B} \\ 
-dA & dA & \Sinv{B} & \Sinv{B} \\ 
\SinvT{B} & \SinvT{B} & dA & -dA \\
\SinvT{B} & \SinvT{B} & -dA & dA \\
\end{pmatrix}
\]
and showed that $V$ is a spin model if and only if $(A,B)$ 
is an invertible Jones pair.
This extends Nomura's result in \cite{N_Twist} which covers
only the invertible Jones pairs $(A,B)$ where both $A$ and $B$ are
symmetric.
We give below a different proof for Yamada's result.

First, it is straightforward to check that $V$ is also a type-II matrix.
Let $\BM$ be the Bose-Mesner algebra of order $4n$ defined
in the previous section.

\begin{theorem}
\label{thm_V_in_B}
If $(A,B)$ is an invertible Jones pair and $A$ is symmetric, then $V$
belongs to $\BM$.
\end{theorem}
{\sl Proof.}
Let $H=dA$.
By Equation~(\ref{cond1}), the matrix $K$ paired with $H$ satisfies
\[
\T{A}{B^T}(K^T)^T = \T{A}{B}(dA)=dB.
\]
Since $\T{A}{B^T}$ is an isomorphism and $\T{A}{B^T}(dA)^T=dB$, 
we conclude that $K=dA^T=dA$.
Hence $V$ is equal to $\M(\0,\inv{A},dA)$ and it belongs to $\BM$.\qed

Assume $(A,B)$ is an invertible Jones pair and $A$ is an $\nbyn$
symmetric matrix.  We use the next four lemmas to show that $\BM
\subseteq \Nom{V}$.  If $M=\M(F,G,H)$ in $\BM$, we want to show that
$\Evector{V}{r}{\Sinv{V}}{s}$ is an eigenvector of $M$ for all
$r,s=1,\ldots,4n$.

In the following, we divide $V$ into sixteen $\nbyn$ blocks.  We use
$\y{\alpha}{\beta}$ to denote $\Evector{V}{r}{\Sinv{V}}{s}$ when
$Ve_r$ is the $i$-th column of the $\alpha$-th block and $\Sinv{V}
e_s$ is the $j$-th column of the $\beta$-th block.  We display the
vectors $\y{\alpha}{\beta}$ to make checking the computation easier.

\[
\y{1}{1} = \y{2}{2} = \Evone, \quad
\y{1}{2} = \y{2}{1} = \Evtwo,
\]

\[
\y{3}{3} = \y{4}{4} = \Evthree, \quad
\y{3}{4} = \y{4}{3} = \Evfour,
\]

\[
\y{1}{3} = -\y{2}{4} = \Evfive,
\]

\[
\y{1}{4} = -\y{2}{3} = \Evsix,
\]

\[
\y{3}{1} = -\y{4}{2} = \Evseven,
\]
and
\[
\y{4}{1} = -\y{3}{2} = \Eveight.
\]

\begin{lemma}
\label{lem_calB_in_NV_1}
Let $M=\M(F,G,H)$ be in $\BM$.
Then for $\range{i,j}{1}{n}$,
$\y{1}{1}$, $\y{1}{2}$, $\y{2}{1}$ and $\y{2}{2}$ are
eigenvectors of $M$.
\end{lemma}
{\sl Proof.}
Note that $M \y{1}{1}$ equals
\[
2
\left(
\begin{aligned}
\Theta_A(F) (\Evector{A}{i}{\Sinv{A}}{j})  &+
\T{A}{B}(G) (\Evector{\SinvT{B}}{i}{B^T}{j}) \\
\Theta_A(F) (\Evector{A}{i}{\Sinv{A}}{j})  &+  
\T{A}{B}(G) (\Evector{\SinvT{B}}{i}{B^T}{j}) \\
\Theta_{\Sinv{B}}(F) (\Evector{\SinvT{B}}{i}{B^T}{j})  &+
\T{A}{B}(G^T)^T   (\Evector{A}{i}{\Sinv{A}}{j})\\
\Theta_{\Sinv{B}}(F)  (\Evector{\SinvT{B}}{i}{B^T}{j}) &+
\T{A}{B}(G^T)^T    (\Evector{A}{i}{\Sinv{A}}{j})\\
\end{aligned}
\right)
\]
which in turn equals
\[
2
\left(
\begin{aligned}
\Theta_A\left(\Theta_A(F)\right)_{i,j}  (\Evector{A}{i}{\Sinv{A}}{j}) &+
\T{A}{B}(G) (\Evector{\SinvT{B}}{i}{B^T}{j}) \\
\Theta_A\left(\Theta_A(F)\right)_{i,j}  (\Evector{A}{i}{\Sinv{A}}{j}) &+
\T{A}{B}(G) (\Evector{\SinvT{B}}{i}{B^T}{j}) \\
\Theta_{\SinvT{B}}\left( \Theta_{\Sinv{B}}(F)\right)_{i,j}
 (\Evector{\SinvT{B}}{i}{B^T}{j}) &+
\T{A}{B}(G^T)^T (\Evector{A}{i}{\Sinv{A}}{j})\\
\Theta_{\SinvT{B}}\left( \Theta_{\Sinv{B}}(F)\right)_{i,j}
 (\Evector{\SinvT{B}}{i}{B^T}{j}) &+
\T{A}{B}(G^T)^T (\Evector{A}{i}{\Sinv{A}}{j})\\
\end{aligned}
\right).
\]
Now, we show that $\y{1}{1}$ is an eigenvector of $M$ and 
compute the corresponding eigenvalue, which is the $ij$-th entry
in the $(1,1)$-block of $\Theta_V(M)$.
Since $A$ is symmetric, it follows from Equation~(\ref{eqn_Duality_A}) that
\begin{eqnarray}
\nonumber
\Theta_A\left(\Theta_A(F)\right)&=&
\Theta_{A^T}\left(\Theta_A(F)\right)\\
\nonumber
&=& n F^T\\
\label{eqn_calB_in_NV_1_1}
&=& \Theta_{\SinvT{B}} \left( \Theta_{\Sinv{B}}(F) \right).
\end{eqnarray}
Moreover, applying Theorem~\ref{thm_T_AB_XDX}~(\ref{eqn_T_AB_XDX_e}) with
$R$ equal to $G$, we have
\[
\XDX{\T{A}{B}(G)}{B^T}{\SinvT{B}} =
\DXD{\SinvT{A}}{A^T}{nG^T}.
\]
Since $A$ is symmetric, the above equation is equivalent to
\begin{equation}
\label{eqn_calB_in_NV_1_2}
\T{A}{B}(G)(\Evector{\SinvT{B}}{i}{B^T}{j})
=nG^T_{i,j}(\Evector{A}{i}{\Sinv{A}}{j})
\end{equation}
for $\range{i,j}{1}{n}$.
Similarly, applying
Theorem~\ref{thm_T_AB_XDX}~(\ref{eqn_T_AB_XDX_d}) with
$R$ equals to $G^T$ gives
\[
\XDX{\T{A}{B}(G^T)^T}{\Sinv{A}}{A}=\DXD{B^T}{\SinvT{B}}{nG^T},
\]
which implies
\begin{equation}
\label{eqn_calB_in_NV_1_3}
\T{A}{B}(G^T)^T (\Evector{A}{i}{\Sinv{A}}{j}) =
nG^T_{i,j}(\Evector{\SinvT{B}}{i}{B^T}{j})
\end{equation}
for $\range{i,j}{1}{n}$.
From Equations~(\ref{eqn_calB_in_NV_1_1}), (\ref{eqn_calB_in_NV_1_2})
and (\ref{eqn_calB_in_NV_1_3}), we see that
\begin{eqnarray*}
M \y{1}{1} 
&=& 
\left(
\begin{aligned}
2nF_{j,i}  (\Evector{A}{i}{\Sinv{A}}{j}) &+ 
2nG_{j,i}  (\Evector{A}{i}{\Sinv{A}}{j}) \\
2nF_{j,i}  (\Evector{A}{i}{\Sinv{A}}{j}) &+ 
2nG_{j,i}  (\Evector{A}{i}{\Sinv{A}}{j}) \\
2nF_{j,i}  (\Evector{\SinvT{B}}{i}{B^T}{j}) &+
2nG_{j,i}  (\Evector{\SinvT{B}}{i}{B^T}{j}) \\
2nF_{j,i}  (\Evector{\SinvT{B}}{i}{B^T}{j}) &+
2nG_{j,i}  (\Evector{\SinvT{B}}{i}{B^T}{j}) \\
\end{aligned}
\right)
\\
&=& 2n(F_{j,i} + G_{j,i}) \y{1}{1},
\end{eqnarray*}
and the $(1,1)$-block of $\Theta_V(M)$ is equal to $2n(F^T+G^T)$.
Since $\y{1}{1}=\y{2}{2}$, the $(2,2)$-block of $\Theta_V(M)$ 
is also $2n(F^T+G^T)$.
For $(\alpha,\beta) \in \{(1,2),(2,1)\}$,
\[
M\y{\alpha}{\beta} = 
2
\left(
\begin{aligned}
-\Theta_A(F) (\Evector{A}{i}{\Sinv{A}}{j})  & 
+\T{A}{B}(G)  (\Evector{\SinvT{B}}{i}{B^T}{j}) \\
-\Theta_A(F) (\Evector{A}{i}{\Sinv{A}}{j})  &
+\T{A}{B}(G)  (\Evector{\SinvT{B}}{i}{B^T}{j}) \\
\Theta_{\Sinv{B}}(F) (\Evector{\SinvT{B}}{i}{B^T}{j})  &-
\T{A}{B}(G^T)^T  (\Evector{A}{i}{\Sinv{A}}{j})\\
\Theta_{\Sinv{B}}(F)  (\Evector{\SinvT{B}}{i}{B^T}{j}) &-
\T{A}{B}(G^T)^T   (\Evector{A}{i}{\Sinv{A}}{j})\\
\end{aligned}
\right).
\]
Using the above argument, the $(1,2)$- and $(2,1)$-blocks of
$\Theta_V(M)$ are equal to $2n(F^T-G^T)$.\qed

\begin{lemma}
\label{lem_calB_in_NV_2}
Let $M=\M(F,G,H)$ be in $\BM$.  Then for $\range{i,j}{1}{n}$,
$\y{3}{3}$, $\y{3}{4}$, $\y{4}{3}$ and $\y{4}{4}$ are eigenvectors of
$M$.
\end{lemma}
{\sl Proof.}
We have
$M \y{3}{3}$ equals
\[
2
\left(
\begin{aligned}
\Theta_A(F) (\Evector{\Sinv{B}}{i}{B}{j}) &+
\T{A}{B}(G) (\Evector{A}{i}{\Sinv{A}}{j}) \\
\Theta_A(F) (\Evector{\Sinv{B}}{i}{B}{j}) &+
\T{A}{B}(G) (\Evector{A}{i}{\Sinv{A}}{j}) \\
\Theta_{\Sinv{B}}(F) (\Evector{A}{i}{\Sinv{A}}{j}) &+
\T{A}{B}(G^T)^T (\Evector{\Sinv{B}}{i}{B}{j})\\
\Theta_{\Sinv{B}}(F) (\Evector{A}{i}{\Sinv{A}}{j}) &+
\T{A}{B}(G^T)^T (\Evector{\Sinv{B}}{i}{B}{j})\\
\end{aligned}
\right)
\]
which is equal to
\[
2
\left(
\begin{aligned}
\Theta_{\Sinv{B}}\left(\Theta_A(F)\right)_{i,j} 
(\Evector{\Sinv{B}}{i}{B}{j}) &+
\T{A}{B}(G) (\Evector{A}{i}{\Sinv{A}}{j}) \\
\Theta_{\Sinv{B}}\left(\Theta_A(F)\right)_{i,j} 
(\Evector{\Sinv{B}}{i}{B}{j}) &+
\T{A}{B}(G) (\Evector{A}{i}{\Sinv{A}}{j}) \\
\Theta_A\left( \Theta_{\Sinv{B}}(F) \right)_{i,j}
(\Evector{A}{i}{\Sinv{A}}{j}) &+
\T{A}{B}(G^T)^T (\Evector{\Sinv{B}}{i}{B}{j})\\
\Theta_A\left( \Theta_{\Sinv{B}}(F) \right)_{i,j}
(\Evector{A}{i}{\Sinv{A}}{j}) &+
\T{A}{B}(G^T)^T (\Evector{\Sinv{B}}{i}{B}{j})\\
\end{aligned}
\right).
\]
We now show that $\y{3}{3}$ is an eigenvector of $M$, and 
compute the corresponding eigenvalue which is the $ij$-th entry
in the $(3,3)$-block of $\Theta_V(M)$.
By Corollary~\ref{cor_TA_TB},
\begin{eqnarray}
\nonumber
\Theta_{\Sinv{B}}\left(\Theta_A(F)\right)
&=&\Theta_{B}\left(\Theta_A(F)\right)^T\\
\nonumber
&=&\inv{B} \Theta_A\left(\Theta_A(F)\right) B\\
\label{eqn_calB_in_NV_2_1}
&=& n \inv{B} F^T B.
\end{eqnarray}
Applying Corollary~\ref{cor_TA_TB} to the Jones pair $(A, B^T)$,
\begin{eqnarray*}
\Theta_A\left( \Theta_{\Sinv{B}}(F) \right)
&=& B^T \Theta_{B^T}\left( \Theta_{\Sinv{B}}(F) \right)^T \invT{B}\\
&=& B^T \Theta_{\SinvT{B}}\left( \Theta_{\Sinv{B}}(F) \right) \invT{B}\\
&=& n B^T F^T \invT{B}\\
&=& n \inv{B} (BB^T) F^T \invT{B}.
\end{eqnarray*}
Since $B\in \dualN{A}{B}$, it follows from
Theorem~\ref{thm_NA_NAB} that $BB^T \in \Nom{A}$.
Now $F^T$ belongs to $\Nom{A}$, the commutativity of $\Nom{A}$ implies
\begin{eqnarray}
\nonumber
\Theta_A\left( \Theta_{\Sinv{B}}(F) \right)
&=& n \inv{B} F^T (BB^T) \invT{B}\\
\label{eqn_calB_in_NV_2_2}
&=& n \inv{B} F^T B.
\end{eqnarray}

From Theorem~\ref{thm_NA_NAB}, there exists $G_1^T \in \N{A}{B^T}$ such that
\[
\T{A}{B}(G)=\T{A}{B^T}(G_1^T)^T.
\]
Hence $G_1$ is paired with $G$.
Applying Theorem~\ref{thm_T_AB_XDX}~(\ref{eqn_T_AB_XDX_d}) with $R$ equals
to $G_1^T$ in $\N{A}{B^T}$ yields
\[
\XDX{\T{A}{B^T}(G_1^T)^T}{\SinvT{A}}{A^T}=\DXD{B}{\Sinv{B}}{nG_1^T},
\]
which is equivalent to
\[
\XDX{\T{A}{B}(G)}{\Sinv{A}}{A}=\DXD{B}{\Sinv{B}}{nG_1^T}.
\]
Consequently
\begin{equation}
\label{eqn_calB_in_NV_2_3}
\T{A}{B}(G)(\Evector{A}{i}{\Sinv{A}}{j}) =
n(G_1^T)_{i,j} (\Evector{\Sinv{B}}{i}{B}{j}).
\end{equation}
By Lemma~\ref{lem_cond2}, $G_1^T$ is also paired with $G^T$.
Applying Theorem~\ref{thm_T_AB_XDX}~(\ref{eqn_T_AB_XDX_e}) to $R=G_1$
in $\N{A}{B^T}$ gives
\[
\XDX{\T{A}{B^T}(G_1)}{B}{\Sinv{B}}=\DXD{\SinvT{A}}{A^T}{nG_1^T},
\]
which is equivalent to
\[
\XDX{\T{A}{B}(G^T)^T}{B}{\Sinv{B}}=\DXD{\Sinv{A}}{A}{nG_1^T}
\]
and
\begin{equation}
\label{eqn_calB_in_NV_2_4}
\T{A}{B}(G^T)^T (\Evector{\Sinv{B}}{i}{B}{j}) = 
n(G_1^T)_{i,j} (\Evector{A}{i}{\Sinv{A}}{j}).
\end{equation}
It follows from Equations~(\ref{eqn_calB_in_NV_2_1}),
(\ref{eqn_calB_in_NV_2_2}), (\ref{eqn_calB_in_NV_2_3})
and (\ref{eqn_calB_in_NV_2_4}) that
\begin{eqnarray*}
M \y{3}{3} 
&=& 2\left(
\begin{aligned}
n(\inv{B}F^TB)_{i,j}  (\Evector{\Sinv{B}}{i}{B}{j}) &+
n(G_1^T)_{i,j}  (\Evector{\Sinv{B}}{i}{B}{j})\\
n(\inv{B}F^TB)_{i,j}  (\Evector{\Sinv{B}}{i}{B}{j}) &+
n(G_1^T)_{i,j}  (\Evector{\Sinv{B}}{i}{B}{j})\\
n(\inv{B}F^TB)_{i,j}  (\Evector{A}{i}{\Sinv{A}}{j}) &+
n(G_1^T)_{i,j}  (\Evector{A}{i}{\Sinv{A}}{j}) \\
n(\inv{B}F^TB)_{i,j}  (\Evector{A}{i}{\Sinv{A}}{j}) &+
n(G_1^T)_{i,j}  (\Evector{A}{i}{\Sinv{A}}{j}) \\
\end{aligned}
\right)
\\
&=& 2n(\inv{B}F^TB+G_1^T)_{i,j} \y{3}{3}.
\end{eqnarray*}
Note that $\y{3}{3}=\y{4}{4}$.
Hence the $(3,3)$- and $(4,4)$-blocks of $\Theta_V(M)$ are equal to
$2n(\inv{B} F^T B + G_1^T)$.
It is easy to see from the block structure of $\y{3}{4}$ and $\y{4}{3}$
that the $(3,4)$- and $(4,3)$-blocks of $\Theta_V(M)$ are equal to
$2n(\inv{B} F^T B -G_1^T)$.\qed

\begin{lemma}
\label{lem_calB_in_NV_3}
Let $M=\M(F,G,H)$ be in $\BM$.  Then for $\range{i,j}{1}{n}$,
$\y{1}{3}$, $\y{1}{4}$, $\y{2}{3}$ and $\y{2}{4}$ are eigenvectors of
$M$.
\end{lemma}
{\sl Proof.}
We have
\begin{eqnarray*}
M \y{1}{3} &=& 2
\left(
\begin{aligned}
d H & (\Evector{A}{i}{B}{j}) \\
-d H & (\Evector{A}{i}{B}{j}) \\
d^{-1} K& (\Evector{\SinvT{B}}{i}{\Sinv{A}}{j}) \\
-d^{-1} K& (\Evector{\SinvT{B}}{i}{\Sinv{A}}{j}) \\
\end{aligned}
\right)\\
&=& 2
\left(
\begin{aligned}
\T{A}{B}(H)_{i,j} & (d\Evector{A}{i}{B}{j}) \\
\T{A}{B}(H)_{i,j} & (-d\Evector{A}{i}{B}{j}) \\
\T{\SinvT{B}}{\Sinv{A}}(K)_{i,j}&(d^{-1}\Evector{\SinvT{B}}{i}{\Sinv{A}}{j}) \\
\T{\SinvT{B}}{\Sinv{A}}(K)_{i,j}&(-d^{-1}\Evector{\SinvT{B}}{i}{\Sinv{A}}{j}) \\
\end{aligned}
\right).
\end{eqnarray*}
By Corollary~\ref{cor_Sinv_AB}, 
\[
\T{\SinvT{B}}{\Sinv{A}}(K) = \T{B^T}{A}(K^T) = \T{A}{B^T}(K^T)^T.
\]
Since $K$ is paired with $H$, by Equation~(\ref{cond1}), the $(1,3)$-block
of $\Theta_M(V)$ is 2$\T{A}{B}(H)$.
Similarly, the $(2,4)$-, $(1,4)$-, $(2,3)$-blocks 
of $\Theta_M(V)$ are equal to 2$\T{A}{B}(H)$.\qed

\begin{lemma}
\label{lem_calB_in_NV_4}
Let $M=\M(F,G,H)$ be in $\BM$.
Then for $\range{i,j}{1}{n}$,
$\y{3}{1}$, $\y{3}{2}$, $\y{4}{1}$ and $\y{4}{2}$ are
eigenvectors of $M$.
\end{lemma}
{\sl Proof.}
We have
\begin{eqnarray*}
M \y{3}{1} &=& 2
\left(
\begin{aligned}
d^{-1} H & (\Evector{\Sinv{B}}{i}{\Sinv{A}}{j}) \\
-d^{-1} H & (\Evector{\Sinv{B}}{i}{\Sinv{A}}{j}) \\
d K& (\Evector{A}{i}{B^T}{j}) \\
-d K& (\Evector{A}{i}{B^T}{j}) \\
\end{aligned}
\right)\\
&=& 2
\left(
\begin{aligned}
\T{\Sinv{B}}{\Sinv{A}}(H)_{i,j} 
& (d^{-1}  \Evector{\Sinv{B}}{i}{\Sinv{A}}{j}) \\
\T{\Sinv{B}}{\Sinv{A}}(H)_{i,j} 
& (-d^{-1}  \Evector{\Sinv{B}}{i}{\Sinv{A}}{j}) \\
\T{A}{B^T}(K)_{i,j} & (d \Evector{A}{i}{B^T}{j}) \\
\T{A}{B^T}(K)_{i,j} & (-d \Evector{A}{i}{B^T}{j}) \\
\end{aligned}
\right)\\
\end{eqnarray*}
By Corollary~\ref{cor_Sinv_AB}, we have
\begin{eqnarray*}
\T{\Sinv{B}}{\Sinv{A}}(H) &=& \T{B}{A}(H^T)\\
&=&\T{A}{B}(H^T)^T\\
&=&\T{A}{B^T}(K),
\end{eqnarray*}
and the last equality follows from the fact that $K^T$ is paired with 
$H^T$.
Therefore the $(3,1)$-block of $\Theta_V(M)$ is equal to $2\T{A}{B}(H^T)^T$.
Similarly, the $(4,2)$-, $(4,1)$- and $(3,2)$-blocks 
are equal to $2\T{A}{B}(H^T)^T$.\qed

\begin{theorem}
\label{thm_calB_in_NV}
If $(A,B)$ is an invertible Jones pair 
and $A$ is symmetric, 
then $\BM$ is a subscheme of $\Nom{V}$.
\end{theorem}
{\sl Proof.}
For any $M \in \BM$, we have shown 
in Lemmas~\ref{lem_calB_in_NV_1} to ~\ref{lem_calB_in_NV_4}
that $\y{\alpha}{\beta}$ is an eigenvector of $M$
for all $\alpha, \beta \in \{1,2,3,4\}$ and $i,j \in \{1,\ldots,n\}$.
Thus $M \in \Nom{V}$ and $\BM \subseteq \Nom{V}$.\qed

\begin{corollary}
\label{cor_TV_calB}
The Bose-Mesner algebra $\BM$ is formally self-dual with
duality map $\Theta_V$.
\end{corollary}
{\sl Proof.}
We see from the proof of Lemmas~\ref{lem_calB_in_NV_1} to 
\ref{lem_calB_in_NV_4} that for $M=\M(F,G,H)$ in $\BM$,
$\Theta_V(M)$ equals 
\[
2
\begin{pmatrix}
nF^T+nG^T & nF^T-nG^T & \T{A}{B}(H) & \T{A}{B}(H)\\
nF^T-nG^T & nF^T+nG^T & \T{A}{B}(H) & \T{A}{B}(H)\\
\T{A}{B}(H^T)^T & \T{A}{B}(H^T)^T 
& n \inv{B}F^TB+nG_1^T & n \inv{B}F^TB-nG_1^T \\ 
\T{A}{B}(H^T)^T & \T{A}{B}(H^T)^T 
& n \inv{B}F^TB-nG_1^T & n \inv{B}F^TB+nG_1^T \\ 
\end{pmatrix},
\]
where $G_1$ is paired with $G$, that is
\[
\T{A}{B}(G)=\T{A}{B^T}(G_1^T)^T.
\]
Since $nF^T \in \Nom{A}$ and $\Nom{A}=\dualNom{A}$, there exists a matrix 
$\hat{F} \in \Nom{A}$ such that $nF^T = \Theta_A(\hat{F})$.
By Corollary~\ref{cor_TA_TB}, we have
\[
\inv{B}nF^TB=\inv{B} \Theta_A(\hat{F}) B = \Theta_B(\hat{F})^T=
\Theta_{\Sinv{B}}(\hat{F}).
\]
By Corollary~\ref{cor_NAB}, we have $G^T \in \N{A}{B}$.  
It follows from Lemma~\ref{lem_cond2} that $G_1^T$ is also paired with $G^T$,
whence we have
\[
\Theta_V(M)=\M(2\hat{F},2H,2nG^T)
\] 
belongs to $\BM$.
Moreover, the map $\Theta_V$ restricted to $\BM$ is a duality map of
$\BM$.\qed

We are ready to prove Yamada's result.
\begin{theorem}[\cite{Yamada2}, Theorem~1]
\label{thm_V_Spin}
Let $A$ be a symmetric $\nbyn$ matrix.
Then $(A,B)$ is an invertible Jones pair
if and only if $V$ is a spin model.
\end{theorem}
{\sl Proof.}
Suppose $(A,B)$ is an invertible Jones pair.
By Theorems~\ref{thm_V_in_B} and \ref{thm_calB_in_NV}, 
the matrix $V$ is equal to $\M(\0,\inv{A},dA)$ and 
hence it belongs to $\Nom{V}$.
By Corollary~\ref{cor_TV_calB}, 
\[
\Theta_V(V) = \M(\0,2dA,2n\inv{A}).
\]
If $K$ is paired with $H=2n\inv{A}$, then
\[
\T{A}{B^T}(K^T)^T=2n\T{A}{B}(\inv{A})=2n\Sinv{B},
\]
which implies $K=2n\inv{A}=2\Sinv{A}$.
Therefore
\begin{eqnarray*}
\Theta_V(V) &=& 2
\begin{pmatrix}
\Sinv{A} & -\Sinv{A} & dB & dB\\
-\Sinv{A} & \Sinv{A} & dB & dB\\
dB^T & dB^T & \Sinv{A} & -\Sinv{A} \\
dB^T & dB^T & -\Sinv{A} & \Sinv{A} \\
\end{pmatrix}\\
&=& 2d \Sinv{V}.
\end{eqnarray*}
By Theorem~\ref{thm_T_AB_XDX}, we have
\[
\XDX{V}{\Sinv{V}}{V} = \DXD{\Sinv{V}}{V}{2d \Sinv{V}}.
\]
Since $V$ is symmetric, 
we conclude that $(\frac{1}{2d}V,\Sinv{V})$ is an invertible Jones pair,
which is equivalent to saying $V$ is a spin model.

Conversely, let $V$ be a spin model, 
or equivalently, let $(\frac{1}{2d}V,\Sinv{V})$ be an invertible Jones pair. 
Since the $(1,3)$-block of $\Sinv{V}$ is equal to $B$, we have
\[
V \y{1}{3}=2d B_{i,j}\y{1}{3}.
\]
This equation implies that
\[
A (\Evector{A}{i}{B}{j}) = B_{i,j} (\Evector{A}{i}{B}{j})
\quad \text{for all } \range{i,j}{1}{n}.
\]
By Theorem~\ref{thm_T_AB_XDX}, we have
\[
\XDX{A}{B}{A}=\DXD{B}{A}{B}.
\]
Similarly, the $(3,1)$-block of $\Sinv{V}$ is equal to $B^T$, we get
\[
V \y{3}{1}=2d (B^T)_{i,j}\y{3}{1},
\]
which implies 
\[
A (\Evector{A}{i}{B^T}{j}) = B^T_{i,j} (\Evector{A}{i}{B^T}{j})
\quad \text{for all } \range{i,j}{1}{n},
\]
and
\[
\XDX{A}{B^T}{A}=\DXD{B^T}{A}{B^T}.
\]
Thus $(A,B)$ is an invertible Jones pair.\qed

It follows from Theorem~\ref{thm_Spin} and the above theorem that
the Bose-Mesner algebra $\Nom{V}$ is formally self-dual and 
$\Theta_V$ is a duality map of $\Nom{V}$.

Given any invertible Jones pair $(C,B)$, it is easy find an odd-gauge
equivalent invertible Jones pair $(A,B)$ in which $A$ is symmetric,
see Section~8 of \cite{CGM}.  By the above theorem, we can always
construct a symmetric spin model $V$ from every invertible Jones pair,
or equivalently, every four-weight spin model.

\section{Subschemes and Induced Schemes}
\label{section_Schemes}

Suppose $A$ and $B$ are $\nbyn$ type-II matrices.
It is easy to verify that the $2n \times 2n$ matrix
\[
W= 
\begin{pmatrix}
A & \Sinv{B} \\ -A & \Sinv{B}\\
\end{pmatrix}
\]
is also a type-II matrix.  Furthermore, if $(A,B)$ is an invertible
Jones pair and $A$ is symmetric, then we have
\begin{equation}
\label{eqn_NW}
\Nom{W} = 
\left\{ 
\begin{pmatrix}
F+G & F-G\\ F-G & F+G\\
\end{pmatrix}
: F \in \Nom{A}, G\in \N{A}{B}
\right\},
\end{equation}
and
\begin{equation}
\label{eqn_NWT}
\Nom{W^T} = 
\left\{ 
\begin{pmatrix}
\Theta_A(F) & \T{A}{B}(G) \\ \T{\Sinv{B}}{\Sinv{A}}(G) & 
 \Theta_{\Sinv{B}}(F)\\
\end{pmatrix}
: F \in \Nom{A}, G\in \N{A}{B}
\right\}.
\end{equation}
Hence the dimensions of $\Nom{W}$ and $\Nom{W^T}$ equal twice the dimension
of $\Nom{A}$.
For details, please see Section~11 of \cite{CGM}.

Now we have five Bose-Mesner algebras $\Nom{V}$, $\BM$, $\Nom{W}$,
$\Nom{W^T}$ and $\Nom{A}$ associated to each invertible Jones pair $(A,B)$
with $A$ symmetric.
The aim of this section is to show that they satisfy the relations
described in the following diagram.
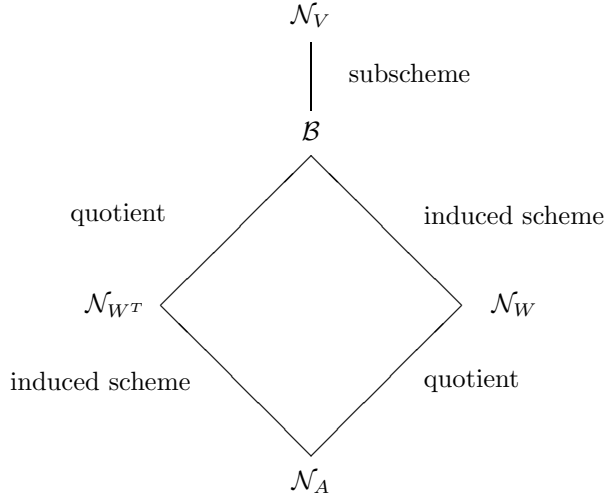
\setlength{\unitlength}{1mm}
\begin{center}
\begin{picture}(85,70)(0,0)
\put(40,10){\line(1,1){20}}
\put(40,10){\line(-1,1){20}}
\put(20,30){\line(1,1){20}}
\put(60,30){\line(-1,1){20}}
\put(40,52){\makebox(0,0)[b]{$\BM$}}
\put(40,5){\makebox(0,0)[b]{$\Nom{A}$}}
\put(10,30){\makebox(0,0)[l]{$\Nom{W^T}$}}
\put(64,30){\makebox(0,0)[l]{$\Nom{W}$}}
\put(55,42){\makebox(0,0)[l]{induced scheme}}
\put(8,42){\makebox(0,0)[l]{quotient}}
\put(55,20){\makebox(0,0)[l]{quotient}}
\put(0,20){\makebox(0,0)[l]{induced scheme}}
\put(40,67){\makebox(0,0)[b]{$\Nom{V}$}}
\put(40,56){\line(0,1){9}}
\put(45,61){\makebox(0,0)[l]{subscheme}}
\end{picture}
\end{center}

Let ${\bf B}$ be a Bose-Mesner algebra on vertex set $\cV$.  Let $Y$ be
a non-empty subset of $\cV$.  For any $|\cV| \times |\cV|$ matrix $M$, we
use $M_Y$ to denote the $|Y| \times |Y|$ matrix obtained from the rows
and the columns of $M$ indexed by the elements in $Y$.  We let the set
\[
{\bf B}_Y:= \{ M_Y : M \in {\bf B}\}.
\]
If ${\bf B}_Y$ is also a Bose-Mesner algebra, we say it is an \textsl{induced
scheme} of ${\bf B}$.  Suppose the vertex sets of $\Nom{A}$, $\Nom{W}$
and $\BM$ are $\{1,\ldots,n\}$, $\{1,\ldots,2n\}$ and
$\{1,\ldots,4n\}$, respectively.  If $Y=\{1,\ldots,n\}$, then it is
obvious from Equation~(\ref{eqn_NWT}) that the set $(\Nom{W^T})_Y$ is
equal to $\dualNom{A}$.  Therefore $\Nom{A}=\dualNom{A}$ is an induced scheme
of $\Nom{W^T}$.  Similarly, let $Y' = \{1,\ldots,2n\}$.  It follows
from Equations~(\ref{eqn_B}) and (\ref{eqn_NW}) that $\BM_{Y'} =
\Nom{W}$.

Let ${\bf B}$ be a Bose-Mesner algebra on vertex set $\cV$.
Let $\pi=(C_1,\ldots, C_r)$ be a partition of $\cV$.
Define the characteristic matrix $S$ of $\pi$ to be the $n \times r$
matrix with
\[
S_{u,k}=
\begin{cases}
1 & \text{if } u \in C_k,\\
0 & \text{otherwise}.
\end{cases}
\]
We say $\pi$ is equitable relative to ${\bf B}$ if and only if for
each matrix $M$ in ${\bf B}$, there is an $r \times r$ matrix $Z_M$
satisfying
\[
MS=SZ_M.
\]
We call the set $\{Z_M : M \in {\bf B}\}$
the \textsl{quotient} of ${\bf B}$ with respect to $\pi$.
For $\range{i}{1}{n}$, let $C_i = \{i,n+i\}$ and let $\pi=(C_1,\ldots,C_n)$.
The characteristic matrix of $\pi$ is
\[
S = 
\begin{pmatrix}
I_n \\ I_n\\
\end{pmatrix}.
\]
Then a matrix 
\[
M=
\begin{pmatrix}
F+R & F-R \\ F-R & F+R\\
\end{pmatrix}
\]
in $\Nom{W}$ satisfies
\[
MS = S (2F).
\]
Thus $Z_M = 2F$.  By Equation~(\ref{eqn_NW}), we see that $F \in
\Nom{A}$ and thus the quotient of $\Nom{W}$ with respect to $\pi$ is
equal to $\Nom{A}$.  Similarly let $C_i=\{i,n+i\}$, for
$i=1,\ldots,n,2n+1,\ldots,3n$.  The characteristic matrix of
$\pi'=(C_1,\ldots,C_n,C_{2n+1},\ldots,C_{3n})$ is
\[
S' = 
\begin{pmatrix}
I_n  & \0 \\ I_n & \0\\
\0 & I_n \\ \0 & I_n \\
\end{pmatrix}.
\]
Then a matrix $\M(F,G,H)$ in $\BM$ satisfies
\begin{eqnarray*}
\M(F,G,H) S' &=& 2
\begin{pmatrix}
\Theta_A(F) & \T{A}{B}(G)\\
\Theta_A(F) & \T{A}{B}(G)\\
\T{A}{B}(G^T)^T & \Theta_{\Sinv{B}}(F)\\
\T{A}{B}(G^T)^T & \Theta_{\Sinv{B}}(F)\\
\end{pmatrix}\\
&=& S' \left(2 
\begin{pmatrix}
\Theta_A(F) & \T{A}{B}(G)\\
\T{A}{B}(G^T)^T & \Theta_{\Sinv{B}}(F)\\
\end{pmatrix}
\right).
\end{eqnarray*}
By Corollary~\ref{cor_Sinv_AB}, we have
\[
\T{\Sinv{B}}{\Sinv{A}}(G) = \T{A}{B}(G^T)^T.
\]
As a result, $Z_{\M(F,G,H)} \in \Nom{W^T}$ and $\Nom{W^T}$ is 
the quotient of $\BM$ with respect to $\pi'$.

In addition, it is straightforward to check that the span of the
following set
\[
\{ \M(F,\0,H) : F \in \Nom{A} \text{ and } H \in \N{A}{B}\}
\cup \{ \M(\0,I_n,\0) \}
\]
is also a Bose-Mesner algebra.  Therefore it is a \textsl{subscheme} of $\BM$
whose dimension equals $2 \dim(\Nom{A}) + 1$.
Similarly,  the span of the set
\[
\left\{
\begin{pmatrix}
\Theta_A(F) & \0 \\ \0 & \Theta_{\Sinv{B}}(F)\\
\end{pmatrix}
: F \in \Nom{A}
\right\}
\cup 
\left\{
\begin{pmatrix}
\0 & J_n \\ J_n & \0 \\
\end{pmatrix}
\right\}
\]
is a subscheme of $\Nom{W^T}$ whose dimension equals $\dim(\Nom{A})+1$.

\section{Comments}
\label{section_Comments}

We now give an explicit description of $\Nom{V}$.
Let ${\cal R}$ be the space consisting matrices
\[
\begin{pmatrix}
\0 & \0 & N & -N\\
\0 & \0 & -N & N\\
N_1 & -N_1 & \0 & \0 \\
-N_1 & N_1 & \0 & \0 \\
\end{pmatrix},
\]
where $N$ and $N_1$ satisfy
\begin{eqnarray*}
\XDX{\inv{A}}{\Sinv{B}}{N} \D_{\Sinv{A}} X_{\inv{B}} &=&  \D_S =
\XDX{B}{A}{N_1} \D_B X_A,\\
\XDX{B^T}{A}{N} \D_{B^T} X_A &=&
\D_{S_1} =
\XDX{\inv{A}}{\SinvT{B}}{N_1} \D_{\Sinv{A}} X_{\invT{B}},
\end{eqnarray*}
for some $\nbyn$ matrices $S$ and $S_1$.  Then $\Nom{V}$ is equal to
the direct sum of $\BM$ and ${\cal R}$, see page~124 of \cite{Chan}.
We see that if $\Nom{A}$ has dimension $r$, then $\dim(\Nom{V})$ equals $3r +
\dim({\cal R})$.  Unfortunately, we do not yet know how to determine
the dimension of ${\cal R}$.  We can only conclude that $3r \leq
\dim(\Nom{V}) \leq 3r+n$.  For example, for each of the three $4
\times 4$ four-weight spin models given in Section~5 of \cite{BB_4wt},
the algebra $\Nom{A}$ has dimension $4$ and $\Nom{V}$ has dimension
$16$.  
The natural problem is to determine the dimension of $\Nom{V}$
for any invertible Jones pair $(A,B)$.  

We get two link invariants from an invertible Jones pair $(A,B)$: 
one from $(A,B)$ and the other from the spin model $V$.  
It is natural to ask how the two invariants are related.
In addition,
it would be very useful to have a procedure that decides whether 
any $4n \times 4n$ spin model is gauge equivalent to a spin model 
that has the same structure as $V$.
Such procedure may lead us to the extraction of  invertible Jones pairs
from the spin models of order divisible by four.


Any new examples of invertible Jones pair will be extremely desirable
since there is a rich family of Bose-Mesner algebras attached.  On the
other hand, we are also interested in any Bose-Mesner algebras that
fit the diagram in Section~\ref{section_Schemes} because they may lead
to the discovery of new invertible Jones pairs, hence possibly new
link invariants.  In particular, we have examined the formally-dual
pair of Bose-Mesner algebras, $\BM_1$ and $\BM_2$, constructed from
the \textsl{Kasami codes} in \cite{Kasami}.  These algebras consist of
$2^{4t+2} \times 2^{4t+2}$ matrices and they have dimension six.  The
Schur-idempotents of, say, $\BM_1$ have valencies
\[
1, 2^{2t+1}-1, 2^{2t+1}-1, 2^{2t+1}-1,
(2^{2t}-1)(2^{2t+1}-1), \text{ and } (2^{2t}-1)(2^{2t+1}-1);
\]
while the valencies of the Schur-idempotents of $\BM_2$ are
\begin{eqnarray*}
1, 2^{2t+1}-1, 
2^{t-1} (2^{t}-1)(2^{2t+1}-1), 2^{t-1} (2^{t}-1)(2^{2t+1}-1),\\ 
2^{t-1} (2^{t}+1)(2^{2t+1}-1), \text{ and } 2^{t-1} (2^{t}+1)(2^{2t+1}-1).
\end{eqnarray*}
We are interested in these algebras because they are the only known
example of a formally-dual pair of Bose-Mesner algebras
that are not translation schemes.
They are candidates for $\Nom{W}$ and $\Nom{W^T}$ in our diagram.

In the following, we use the structure of $\Nom{W^T}$ 
to rule out the possibility that $\BM_1$ and $\BM_2$ fit 
into the diagram in Section~\ref{section_Schemes}.
We see from the previous section that 
\[
\hat{J} = 
\begin{pmatrix}
\0 & J_{2^{4t+1}} \\ J_{2^{4t+1}} & \0\\
\end{pmatrix}
\]
belongs to $\Nom{W^T}$.  Therefore if $\Nom{W^T}$ equals to $\BM_1$, then a
subset of the Schur-idempotents of $\BM_1$ would sum to $\hat{J}$.  In
this case, a subset of the valencies of $\BM_1$ would sum to
$2^{4t+1}$.  However, we can use elementary computation to prove that
it is impossible to find a subset of the numbers in the first list
above to sum to $2^{4t+1}$.  Consequently the algebra $\BM_1$ cannot
be $\Nom{W^T}$.  Similarly, simple computation shows that we cannot
find a subset of valencies of $\BM_2$ to sum to $2^{4t+1}$.  We conclude that
$\BM_2$ cannot be $\Nom{W^T}$.  As a result there does not exist any
invertible Jones pair for which $\{\BM_1,\BM_2\}$ equals
$\{\Nom{W},\Nom{W^T}\}$.

\section*{Acknowledgement}
We thank the referees for their constructive comments and suggestions.

\bibliographystyle{acm}
\bibliography{JP}
\end{document}